\documentclass[11pt,leqno]{article}
\usepackage{amsmath,amssymb,latexsym}
\usepackage{amscd}
\usepackage{mathrsfs}
\usepackage{theorem}
\usepackage{amsfonts}
\usepackage{color}
\newtheorem{proposition}{Proposition}[section]
\newtheorem{example}[proposition]{Example}
\newtheorem{lemma}[proposition]{Lemma}
\newtheorem{corollary}[proposition]{Corollary}
\newtheorem{theorem}[proposition]{Theorem}
\newtheorem{remark}[proposition]{Remark}

\newtheorem{problem}[proposition]{Problem}

\newcommand{\Proof}{\textbf{Proof:} \ }
\newcommand{\qed}{\hspace*{\fill} $\Box $}

\newcommand{\cA}{\mathscr{A}}

\newcommand{\bC}{\mathbb{C}}
\newcommand{\bN}{\mathbb{N}}
\newcommand{\bR}{\mathbb{R}}

\newcommand{\bT}{\mathbb{T}}

\newcommand{\vp}{\varphi}
\newcommand{\eps}{\varepsilon}

\newcommand{\Ar}{\mbox{$\cA (\bR)$}}

\newcommand{\id}{\operatorname{id}\,}
\newcommand{\im}{\operatorname{im}\,}

\newcommand{\nN}{{n\in \bN}}

\newcommand{\cf}{C_\vp}

\title{{\sc A note on the spectrum of composition operators on  spaces of real analytic functions}}
\author{{\sc  Jos\'e Bonet and Pawe{\l} Doma{\'n}ski}}

\date{}

\begin{document}

\maketitle

\vspace{2cm}

{\bf Authors' Addresses:}

\vspace{2mm}

\begin{minipage}{7.5cm}
J.\ Bonet (corresponding author)

Instituto Universitario de Matem\'{a}tica Pura y Aplicada IUMPA

Universitat Polit\`{e}cnica de Val\`{e}ncia

E-46071 Valencia, SPAIN

e-mail: jbonet@mat.upv.es

phone: +34963879497

fax: +34963879497
\end{minipage}
\hspace*{\fill}
\parbox{7.5cm}{P. Doma{\'n}ski

Faculty of Mathematics and Comp. Sci.

A. Mickiewicz University Pozna{\'n}

Umultowska 87

61-614 Pozna{\'n}, POLAND

e-mail: domanski@amu.edu.pl

}

\vspace{2cm}

\begin{abstract} In this paper the spectrum of composition operators on the space of real analytic functions is investigated. In some cases it is completely determined while in some other cases it is only estimated.
\end{abstract}

\footnotetext[1]{


{\rm This paper is accepted for publication in Complex Analysis and Operator Theory.} 
}

\newpage

\section{Introduction}

Let $\vp:\bR\to\bR$ be a non-constant real analytic map and let $\Ar$ be the space of real analytic functions
defined on $\bR$. Each symbol $\vp:\bR\to\bR$ defines a composition operator $\cf: \Ar \to \Ar$ by $\cf(f):= f \circ \vp, f \in \Ar$.
When $\Ar$ is endowed with its natural locally convex topology (see e.g.\ \cite{Dwin}), $\cf$ is a continuous linear operator on $\Ar$. In our article \cite{BDeigen} we studied the eigenvalues and eigenvectors of composition operators $\cf: \Ar \to \Ar$. In this note we complement those results with some examples and remarks concerning the spectrum of $\cf: \Ar \to \Ar$.

Our results in \cite{BDeigen} give  precise information about the
injectivity of the operator $\cf - \lambda I: \Ar \to \Ar$,
$\lambda \in \bC \setminus \{0\}$. In order to investigate the
spectrum $\sigma(\cf)$, we have to study when the operator $\cf -
\lambda I: \Ar \to \Ar$ is surjective. This operator is surjective
if and only if the operator $T_{\lambda}: \Ar \to \Ar$, defined by
$T_{\lambda}f(x):= f(x) - (1/\lambda) f(\vp(x)), x \in \bR, f \in
\Ar,$ is surjective. If we set $g(x,y):=(1/\lambda)y + \gamma(x)$,
for $\gamma \in \Ar$, then $T_{\lambda}: \Ar \to \Ar$ is
surjective if and only if for each $\gamma \in \Ar$ there is $f
\in \Ar$ such that $f(x)=g(x,f(\vp(x)))$, $x \in \bR$. We use
methods of the paper \cite{BT}.

Equations of the form $f(x)=g(x,f(\vp(x))$ with $g: \bR^2 \to \bR$
and $\vp: \bR \to \bR$ real analytic are considered in \cite{BT}
and \cite{Sm} (see also \cite{BL1}, \cite{BL2} or the book
\cite{BTb}). As it is mentioned in the introduction of \cite{BT},
if $u$ is a fixed point of $\vp$ and $\lambda \neq \vp'(u)^n$ for
each $n \in \bN_0$, then $f(x) = (1/\lambda) f(\vp(x)) +
\gamma(x)$ has a formal solution $f(x)= u + \sum_{n=1}^{\infty}
c_n(x-u)^n$. Smajdor \cite{Sm} studies conditions to ensure that
this series representation has a positive radius of convergence
near $u$.

In the rest of the article we denote $\id(x)=x, x \in \bR,$ and,
for a map $\vp:\bR\to\bR$, we write $\vp^{[0]}=\id$ and
$\vp^{[n]}$ for the $n$-times composition of $\vp$, $n \in \bN$.
By $I:\cA(J)\to \cA(J)$ we denote the identity operator.

The composition operators on spaces of real analytic functions have been considered
in several papers, like \cite{Dwin}, \cite{DGL}, \cite{DL}, \cite{DL2},
\cite{DL3}, \cite{BDradyn}, \cite{BDhy}.  For  literature on the
space of real analytic functions see a recent survey \cite{Dwin}.
 For
functional analytic tools see \cite{MV}.

\section{Preliminaries} By interval $J$ we mean also unbounded
ones (i.e., halflines or the whole real line). It is clear that
results on the (point) spectrum of a composition operator on $\Ar$
can be easily transferred to the case of $\cA(J)$ where $J$ is an
open interval in $\bR$. We recall the following two results from
\cite{BDeigen}.

\begin{proposition}\label{spectrum0} Let $\vp:J \to J$ be a real analytic function on an open interval
$J \subset \bR$ and let $\cf: \cA(J) \to \cA(J)$ be the associated
composition operator. Then

(1) (\cite[Proposition 1.1 (4)]{BDeigen}) $0 \in \sigma(\cf)$ if
and only if $\vp:J \to J$ is not a real analytic diffeomorphism.

(2) (\cite[Proposition 2.8 (4)]{BDeigen})  If $u \in J$ is a fixed
point of $\vp$ such that $|\vp'(u)| \neq 0, 1$, then $\vp'(u)^n
\in \sigma(\cf)$ for each $n \in \bN_0$.
\end{proposition}

 We say that some subset $A$ of an open interval $J\subset \bR$ is {\em bounded away from the upper end of} $J$ if there is $\delta\in J$ such that all elements of $A$ are $<\delta$. Analogously we define subsets {\em bounded away form the lower end of} $J$.  We also use the following description of the point spectrum of
$\cf$ from \cite{BDeigen}:

\begin{theorem}\label{eigen}
Let $\vp:J\to J$ be a real analytic map for some open interval
$J\subset \bR$.
\begin{itemize}
\item[(a)] If $\vp$ has no fixed points and the set of its
critical points is bounded away from the upper end of $J$ (in case $\vp>\id$) or from
the lower end of $J$ (in case $\vp<\id$) then $\sigma_p(\cf)=\bC\setminus\{0\}$
and every eigenspace is isomorphic to $\cA(\bT)$.
\item[(b)] If $\vp$ has a fixed point $u$ then:
\begin{itemize}
\item[(b1)] if $\vp^{[2]}$ has exactly one fixed point and either
$0<\vert \vp'(u)\vert<1$ or $1<\vert \vp'(u)\vert$ and $\vp$ has
no critical points then
$$
\sigma_p(\cf)=\{(\vp'(u))^n:n=0,1,\dots\}
$$
and the eigenspaces are all one dimensional;
\item[(b2)] if $\vp^{[2]}=\id\neq \vp$ then
$\sigma_p(\cf)=\{-1,1\}$ and the eigenspaces are isomorphic to
$\cA_+(\bR)$ the space of even real analytic functions;
\item[(b3)] if $\vp=\id$ then $\sigma_p(\cf)=\{1\}$ and the
eigenspace is equal to $\cA(J)$.
\end{itemize}
\item[(c)] In all other cases $\sigma_p(\cf)=\{1\}$ and the
eigenspace consists of constant functions only.
\end{itemize}
\end{theorem}

We will use also the following result of Smajdor \cite[Theorem p.
40]{Sm}:

\begin{theorem}\label{Smaj} Let $\vp$ be a holomorphic function of
one complex variable on a neighbourhood of $u\in \bC$, $\vp(u)=u$,
$\vert \vp'(u)\vert<1$. Let $h$ be a holomorphic  function of two
complex variables on a neighbourhood of $(u,v)\in \bC^2$,
$h(u,v)=v$. If there is a formal solution $f(z)=\sum_n f_n(z-u)^n$
of the equation
\begin{equation}\label{eq1}
f(z)=h(z,f(\vp(z)), \qquad f(u)=v,
\end{equation}
then there is a holomorphic solution of (\ref{eq1}) on some
neighbourhood of $u$.

In case
\begin{equation}\label{eq2}
1-\frac{\partial h}{\partial v}(u,v)\cdot(\vp'(u))^n\neq 0\qquad
\text{for every $n=0,1,\dots$}
\end{equation}
the formal solution is unique and so it is convergent around $u$
and gives a holomorphic solution of (\ref{eq1}) on some
neighbourhood of $u$.
\end{theorem}

Please note that in the proof of Proposition \ref{case2} we give a simple example where for $\vert \vp'(u)\vert=1$ the result above fails.

\section{Spectrum}
The case of $\vp:J\to J$ without fixed points is mostly solved.

\begin{corollary}\label{fpf}
Let $\vp:J\to J$ be real analytic, $J\subset \bR$ an open
interval, and $\vp$ have no fixed points.
\begin{itemize}
\item[(a)] If $\vp$ is a diffeomorphism onto $J$ then
$\sigma(\cf)=\sigma_p(\cf)=\bC\setminus \{0\}$.
\item[(b)] If $\vp$ is not a diffeomorphism onto $J$ but the set
of critical points of $\vp$ is bounded away from the upper end of $J$ (in case
$\vp>\id$) or from lower end of $J$ (in case $\vp<\id$) then
$\sigma(\cf)=\bC$ and $\sigma_p(\cf)=\bC\setminus \{0\}$.
\end{itemize}
\end{corollary}

\Proof Combine Proposition \ref{spectrum0} and Theorem
\ref{eigen}.\qed

\begin{problem}
Describe $\sigma(\cf)$ if $\vp>\id$ has no fixed points but the
set of critical points of $\vp$ is not bounded away from the upper end of the domain.
\end{problem}

Now, we concentrate on the fixed point case. We will need the
following standard extension lemma (this procedure for
diffeomorphisms was  used in \cite{BT}, we write it
precisely in the general case).

\begin{lemma}\label{ext}
Let $\vp:(a,b)\to (a,b)$ be a real analytic map on an open
interval $(a,b)\subset \bR$ and let $\vp(c,d)\subset (c,d)\subset [c,d] \subset
(a,b)$ for some fixed $a,b,c,d$. Assume that for every $x\in
(a,b)$ there is $n\in \bN$ such that $\vp^{[n]}(x)\in (c,d)$
(i.e., $(a,b)$ is the attraction basin for $(c,d)$).

If $\tilde f\in \cA(c,d)$ satisfies
\begin{equation}\label{eq3}
\tilde f(\vp(x))-\lambda \tilde f(x)=\gamma(x)\qquad \text{for every
$x\in (c,d)$},
\end{equation}
for some fixed $\lambda\in \bC\setminus \{0\}$ and $\gamma \in
\cA(a,b)$, then $\tilde f$ extends to $f\in \cA(a,b)$ satisfying
\begin{equation}\label{eq4}
f(\vp(x))-\lambda f(x)=\gamma(x)\qquad \text{for every $x\in
(a,b)$}.
\end{equation}
\end{lemma}

\Proof Define
$$
J_n:=\{x\in (a,b)\mid \vp^{[n]}(x)\in (c,d)\}.
$$
Clearly, $(J_n)$ is an increasing open exhaustion of $(a,b)$. Take
any compact increasing exhaustion $([a_n,b_n])_\nN$        with
$$
(c,d)\subset [a_n,b_n]\subset (a,b)\quad \text{and}\quad
\bigcup_\nN [a_n,b_n]=(a,b).
$$
Without loss of generality we may assume that $J_n\supset
[a_n,b_n]$.

We define inductively:
$$
f_0:=\tilde f, \qquad
f_n(x):=\frac{1}{\lambda}\left(\gamma(x)+f_{n-1}(\vp(x))\right)\qquad
\text{for $x\in J_n$}.
$$
By (\ref{eq3}) it is easy to observe that $f_n=\tilde f$ on
$(c,d)$ and thus $f_n$ extends $\tilde f$ on $[a_n,b_n]$. Now,
defining $f(x):=f_n(x)$ for $x\in [a_n,b_n]$ we obtain a real
analytic function $f\in \cA(a,b)$. Since
$$
f(\vp(x))-\lambda f(x)=\gamma(x) \qquad \text{for every $x\in
(c,d)$}
$$
and since both sides of the above equality are analytic on $(a,b)$
so the equality holds everywhere on $(a,b)$. \qed

Note that if $\vp^{[2]}$ has the unique fixed point $u$ then $u$ must be the unique fixed point of $\vp$ as well since otherwise $\vp(u)=w\neq u$ but then $\vp^{[2]}(w)=w$ and $\vp^{[2]}$ would have had two fixed points $u$ and $w$.

\begin{theorem}\label{spectrum1}
Let $\vp:J \to J$ be a real analytic function on an open interval
$J \subset \bR$  and let $\cf: \cA(J) \to
\cA(J)$ the associated composition operator. Suppose that
$\vp^{[2]}$ has a unique fixed point $u$ and $|\vp'(u)| < 1$. Then

(i) If $\lambda \neq 0$ and $\lambda \neq \vp'(u)^n$ for each $n
\in \bN_0$, then $\lambda \notin \sigma(\cf)$, and

(ii) $\ker(\cf - \lambda I)$ is finite dimensional for all
$\lambda \in \bC$.
\end{theorem}

\begin{remark} The case of $\vp$ diffeomorphic onto is proved in
\cite[Th. 4.4]{BT}.
\end{remark}

\Proof (i): Consider $\lambda \neq 0$, $\lambda \neq \vp'(u)^n, n
\in \bN_0$. By Theorem \ref{eigen}, $\cf-\lambda I$ is injective,
we will show that it is surjective.

It suffices to solve, for every $\gamma \in \cA(J)$, the equation
$$f(x) - (1/\lambda) f(\vp(x)) = \gamma(x),$$
which is equivalent to
$$f(x) = g(x,f(\vp(x)),$$
with $g(x,y)=\gamma(x) + (1/\lambda) y$. In order to apply Theorem
\ref{Smaj} to the equation above, we define $v:=\lambda
\gamma(u)/(\lambda -1)$. Clearly, $g(u,v)=v$. Since
$$
\frac{\partial g}{\partial v}(u,v)=\frac{1}{\lambda},
$$
the condition (\ref{eq2}) is satisfied.

 Therefore all the
assumptions of Theorem \ref{Smaj} are satisfied and the unique
formal solution $f_0$ is real analytic on a neighbourhood of $u$.
Hence there is $r>0$ such that $[u-r,u+r] \subset J$ and $f_0 \in
\cA(u-r,u+r)$ satisfies
$$
f_0(x)-1/\lambda f(\vp(x))=\gamma(x)\qquad \text{for every $x\in
(u-r,u+r)$}.
$$

Since $|\vp'(u)| < 1$, we may assume without loss generality that
$\vp(]u-r,u+r[) \subset ]u-r,u+r[$. In the proof of \cite[Theorem
2.6]{BDeigen} it is proved that if $\vp:J\to J$, $J\subset \bR$ an
open interval, has a fixed point $u$, $\vert \vp'(u)\vert<1$, and
$\vp^{[2]}$ has exactly one fixed point, then for every $x\in J$
holds $\vp^{[n]}(x)\to u$ as $n\to \infty$. Thus $(c,d)=(u-r,u+r)$
satisfies the assumptions of Lemma \ref{ext}. This completes the
proof.

(ii): This follows from results of \cite[Section 2]{BDeigen} transferring
them to an arbitrary interval instead of the whole line. For readers
convenience we sketch the proof. By assumption $\vp^{[2]}$ has a unique
fixed point in $J$ and we can aply \cite[Theorem 2.6]{BDeigen}
to conclude that $\vp'(u)^n$ is an eigenvalue of $\cf$
for each $n \in \bN_0$. By \cite[Theorem 2.9 (c) and (d)]{BDeigen},
$\ker(\cf - \lambda I)$ is one dimensional for $\lambda=0$ (if
$\vp'(u)=0$) or $\lambda=\vp'(u)^n$ (if $0<|\vp'(u)|<1)$. We can
apply \cite[Propositions 1.1 and 2.4]{BDeigen} to
conclude that $\ker(\cf - \lambda I) = \{ 0 \}$ for $\lambda \neq
\vp'(u)^n, n \in \bN_0$. \qed

\begin{corollary}\label{cor1spectrum1}
If $\vp:J \to J$ is a real analytic  on an open interval $J
\subset \bR$,   such that $\vp^{[2]}$ has the unique
fixed point $u$ satisfying $|\vp'(u)| < 1$, then
$$
\sigma(\cf)=\begin{cases}\{ \vp'(u)^n \ | \ n \in \bN_0 \} &
\text{if $\vp$ is
a diffeomorphism},\\
\{ \vp'(u)^n \ | \ n \in \bN_0 \}\cup \{0\} & \text{otherwise}
\end{cases}
$$
and
$$
\sigma_p(\cf)=\{ \vp'(u)^n \ | \ n \in \bN_0 \}\setminus \{0\}
$$
Moreover, $\ker(\cf - \lambda I)$ is finite dimensional for all
$\lambda \in \bC$.
\end{corollary}

\Proof The description of $\sigma_p(\cf)$ follows from Theorem
\ref{eigen}. The description of $\sigma(\cf)$ follows
combining of Proposition \ref{spectrum0} and Theorem
\ref{spectrum1}.  \qed

\begin{corollary}\label{cor2spectrum1}
If $\vp:J \to J$ is a real analytic diffeomorphism on an open
interval $J \subset \bR$  such that $\vp^{[2]}$ has a unique fixed
point $u$ satisfying $|\vp'(u)| > 1$, then
$\sigma(\cf)=\sigma_{p}(\cf)=\{ \vp'(u)^n \ | \ n \in \bN_0 \}$
and $\ker(\cf - \lambda I)$ is finite dimensional for all $\lambda
\in \bC$.
\end{corollary}

\Proof  The map $\psi:=\vp^{-1}: J \to J$ is a real analytic
diffeomorphism with a unique fixed point $u \in J$ and $0 <
|\psi'(u)| < 1$. The conclusion follows from Corollary
\ref{cor1spectrum1} and the following observation: for $\mu \neq
0$, $\mu \neq \psi'(u)^n = 1/\vp'(u)^n, n \in \bN_0,$ the operator
$C_{\psi} - \mu I: \cA(J) \to \cA(J)$ is surjective/injective if
and only if $\cf - (1/\mu) I: \cA(J) \to \cA(J)$ is
surjective/injective. \qed

\begin{problem} Describe $\sigma(\cf)$ if $\vp$ has a fixed point
$u$, $\vert \vp'(u)\vert>1$, $\vp^{[2]}$ has a unique fixed point but $\vp$
is not a diffeomorphism.
\end{problem}

By \cite[Theorem 4.5]{BT} we obtain immediately:

\begin{proposition}\label{spectrum2}
Let $\vp: j \to J$, $J\subset \bR$ open interval,  be a real
analytic diffeomorphism with fixed points $x_1 < x_2 <...< x_s$
such that $|\vp'(x_i)| \neq 1$ for $i=1,2,...,s$, $s>1$. For each
$\lambda \neq 0$, the operator $\cf -{\lambda}I: \cA(J) \to
\cA(J)$
 has
closed range, $\ker \cf-{\lambda}I = \{ 0 \}$ if $\lambda \neq 1$,
$ \ker \cf-I $ is finite dimensional and the codimension of $\im
\cf-{\lambda}I$ in $\cA(J)$ is infinite. In particular,
$\sigma_{p}(\cf)= \{ 1 \}$ and $\sigma(\cf) = \bC \setminus \{ 0
\}$.
\end{proposition}

\Proof By Proposition \ref{spectrum0} (1), $0 \notin \sigma(\cf)$
and by Theorem \ref{eigen}, $\sigma_{p}(\cf)= \{ 1 \}$, and
$\ker(\cf - I)$ is one dimensional. The rest of the statement is a
direct consequence of \cite[Theorem 4.5]{BT}. \qed

\begin{example}\rm
(a) If $\vp(x)=e^{\alpha x}, \alpha >1/e$, then
$\sigma_{p}(\cf)=\bC \setminus \{ 0 \}$ and $\sigma(\cf)= \bC$
(use Corollary \ref{fpf}).

(b) If $\vp$ has no fixed points on $\bR$ and $\vp$ has only
finitely many critical points, then $\sigma_{p}(\cf)=\bC\setminus
\{ 0 \}$ and $\sigma(\cf)= \bC$ by Corollary \ref{fpf}.

(c) If $\vp(x) = \alpha \arctan(x), 0<\alpha<1$, then
$\sigma_{p}(\cf)= \{ \alpha^n \ | \ n \in \bN_0 \}$ and
$\sigma(\cf)=\sigma_{p}(\cf) \cup \{ 0 \}$. In this case $\vp$ has
a unique fixed point $u=0$ and $\vp'(u)=\alpha$. We have
$0<|\vp'(u)|<1$  and $\vp$ is not surjective, hence the statement
follows from Corollary \ref{cor1spectrum1}.
\end{example}

Before we  deal with more examples, we present for the reader's
convenience the main result in \cite[Section 3]{BT} in a way
suitable for our purposes explaining details omitted in \cite{BT}.

Let $\vp: \bR \to \bR$ be a real analytic self map, let $\lambda
\in \bC, \lambda \neq 0,$ and consider the operator $T: \Ar \to
\Ar, Tf(x):=f(\vp(x))-\lambda f(x)$. Let $(U_j)_j$ be an open
covering of $\bR$ such that $\vp(U_j) \subset U_j$ for each $j$
and such that each $U_j$ intersects only finitely many other
$U_k$'s. Denote by $T_k: \cA(U_k) \to \cA(U_k)$ and $T_{k,j}:
\cA(U_k \cap U_j) \to \cA(U_k \cap U_j)$ the operator $T$ defined
on the corresponding space. With the notation $g|V$ for the
restriction of the function $g$ to the set $V$, define
$$
\begin{aligned}
\cA(\bR,(U_k),T)&:= \{ \gamma \in \Ar \ : \ \gamma|U_k \in
\im(T_k) \ \  \forall k \},\\
E&:=\{ (f_{k,j})_{k,j} \ : \ f_{k,j}
= - f_{j,k}, \ \ f_{k,j} \in \ker(T_{k,j}) \},\\
F&:=\{(f_k)_k \ : \ f_k \in \ker(T_k) \ \ \forall k \},
\end{aligned}
$$
and
$$
S:F \to E, \ \  S((f_k)_k:= ((f_k - f_j)|U_k \cap U_j)_{k,j}.
$$

\begin{theorem}\label{theoremBT} (\cite[Theorem 3.1]{BT}) There is a linear surjective map
$$\mathcal{F}: \cA(\bR,(U_k),T) \to E/S(F)$$
with kernel equal to $\im(T)$.

In particular, if $S(F) \neq E$, then $T$ is not surjective.
\end{theorem}

\Proof Given $\gamma \in \cA(\bR,(U_k),T)$, for each $k$ there is
$f_k \in \cA(U_k)$ such that $T_k f_k = \gamma|U_k$. Set
$f_{k,j}:=(f_k - f_j)|U_k \cap U_j$. Clearly $T_{k,j} f_{k,j}(x) =
0$ for each $x \in U_k \cap U_j$ and each $k,j$. Hence
$(f_{k,j})_{k,j} \in E$. We define
$\mathcal{F}(\gamma)=(f_{k,j})_{k,j} + S(F)$. To see that
$\mathcal{F}$ is well defined, suppose that $T_k g_k = \gamma|U_k$
with $g_k \in \cA(U_k)$ and set $g_{k,j}:=(g_k - g_j)|U_k \cap
U_j$. We have $(f_{k,j}-g_{k,j})_{k,j} = S((f_k - g_k)_k)$ and
$T_k(f_k - g_k) = 0$ for each $k$. Therefore $(f_{k,j})_{k,j} -
(g_{k,j})_{k,j} \in S(F)$ and $\mathcal{F}$ is indeed well
defined.

We show now that $\mathcal{F}$ is surjective. Fix $(f_{k,j})_{k,j} + S(F) \in E/S(F)$. Our assumptions on the covering
$(U_k)_k$ of $\bR$ permit us to find open sets $V_k$ in the complex plane with $V_k \cap \bR = U_k$ and extensions of
the functions $f_{k,j}$ to holomorphic functions $g_{k,j}$ on $V_k\cap V_j$.

Indeed, let  $f_{k,j}\in \cA(U_k\cap U_j)$ then there is an open
complex neighbourhood $W$ of $U_k\cap U_j$, $W\cap \bR=U_k\cap
U_j$, such that $f_{k,j}$ extends holomorphically as a function
$g_{k,j}$ on $W$. We define:
$$
\tau_j(x):=\begin{cases} \min(d(x,\partial_{\bR}U_j),
d(x,\partial_{\bC} W))  & \text{if $x\in U_k\cap U_j$ and
$d(x,\partial_{\bR}U_j)\leq d(x, \partial_{\bR}U_k)$, $k\neq
j$};\\
d(x,\partial_{\bR} U_j)  & \text{otherwise}.
\end{cases}
$$
We define $\tau_k$ analogously.
Then we define open sets $V_j$ and $V_k$, for instance,
$$
V_j:=\{x+iy\in \bC^d: x\in U_j, \vert y\vert <\tau_j(x)\}.
$$
Then
$$
V_k\cap V_j\subset W, \quad V_k\cap V_j\cap \bR=U_k\cap
U_j,\quad V_j\cap \bR=U_j,\quad V_k\cap\bR=U_k.
$$
If $U_k$ intersects more sets $U_j$ (but finitely many) then we
can take intersections of the obtained sets $V_k$. In that way we
get $g_{k,j}$ holomorphic on $V_k\cap V_j$.

Since for the orthogonal projection $p:\bC=\bR^2\to \bR$ we have $p(V_j)=V_j\cap \bR=U_j$ (see the definition of $V_j$), thus
$$
p(V_j\cap V_k\cap V_l)=U_j\cap U_k\cap U_l
$$
and every connected component of $V_j\cap V_k\cap V_l$ contains a non-empty open part of $U_j\cap U_k\cap U_l$. As easily seen on the intersection $U_j\cap U_k\cap U_l$ it holds $f_{j,k}+f_{k,l}+f_{l,j}=0$ by the very definition of the functions. Thus also extensions satisfy $g_{jk}+g_{kl}+g_{lj}=0$ on $V_j\cap V_k\cap V_l$ for any triple $\{j,k,l\}$.

By the Cartan-Grauert result (see \cite[Prop. 1]{C}) every open
set in $\bR^d$ has a basis of complex neighbourhoods being domains
of holomorphy. So there is a domain of holomorphy $V\subset
\bigcup V_k$, $V\cap J=J$. Thus we can apply the solution to the first
Cousin problem \cite[Th. 5.5.1]{H} for the covering $(V_k\cap V)$ of $V$ and
Cousin data $(g_{j,k})$ to find, after taking restrictions,
$(f_k)_k$ with $f_k \in \mathcal{A}(U_k)$ such that $f_{k,j}=f_k -
f_j$ for each $k,j$. Define $\gamma_k := T_k f_k \in
\mathcal{A}(U_k)$. Since $ f_{k,j} \in \ker(T_{k,j})$,
$\gamma_k(x) = \gamma_j(x)$ for each $x \in U_k \cap U_j$, and we
can find $\gamma \in \Ar$ such that $\gamma|U_k = \gamma_k =T_k
f_k$ for each $k$. This means $\gamma \in \cA(\bR,(U_k),T)$.
Moreover $\mathcal{F}(\gamma)=(f_{k,j})_{k,j} + S(F)$.

It remains to show that $\ker(\mathcal{F}) = \im(T)$. If $\gamma
\in \im(T)$, there is $f \in \Ar$ such that $Tf=\gamma$. Then
$\gamma \in \cA(\bR,(U_k),T)$ since $\gamma|U_k = T_k(f|U_k)$ for
each $k$. Moreover, $\mathcal{F}(\gamma)= 0 \in E/S(F)$ and
$\gamma \in \ker(\mathcal{F})$.

To prove the reverse inclusion, take $\gamma \in \cA(\bR,(U_k),T)$
with $\mathcal{F}(\gamma)= 0 \in E/S(F)$. For each $k$ there is
$f_k \in \cA(U_k)$ such that $T_k f_k = \gamma|U_k$. By assumption
$((f_k - f_j)|U_k \cap U_j)_{k,j} \in S(F)$. Therefore we can find
$(g_k)_k$, $g_k \in \cA(U_k)$, $T_k g_k = 0$ such that $f_k - f_j
= g_k - g_j$ on $U_k \cap U_j$ for each $k,j$. Hence $f_k - g_k =
f_j - g_j$ on $U_k \cap U_j$ for each $k,j$. Accordingly, there is
$f \in \Ar$ such that $f|U_k = f_k - g_k$ for each $k$. If $x \in
U_k$, we get $Tf(x) = T_k f_k(x) -T_k g_k(x)= \gamma(x)$. Since
$(U_k)$ is an open  covering of $\bR$ we get $Tf = \gamma$ and
$\gamma \in \im(T)$. \qed

\begin{proposition}\label{exmonomials}
If $\vp(x)=x^s, s \in \bN$, $s>1$, then for $\cf:\Ar\to\Ar$ we have
that $\sigma_{p}(\cf)= \{ 1 \}$, $\ker(\cf -I)$ is one dimensional
and its elements are  the constants, but $\sigma(\cf)=
\bC$.
\end{proposition}

\Proof The statement on eigenvalues and eigenspaces follows
by Theorem \ref{eigen}.

Now, we consider the spectrum.
First of all, since $\vp$ is not a diffeomorphism, $0 \in
\sigma(\cf)$ by Proposition \ref{spectrum0}. The map $\vp$ has
three fixed points $0,1,-1$ if $s$ is odd and two $0,1$ if $s$ is
even. Anyway, $\vp'(1)=s$ (which coincides with $\vp'(-1)$ if $s$
is odd). By Proposition \ref{spectrum0}, as $s \neq 1$, $s^n \in
\sigma(\cf)$ for each $n \in \bN_0$. It remains to show that
$$
T:=\cf -\lambda I: \Ar \to \Ar
$$
is not surjective for $\lambda \neq 0, \lambda \neq s^n, n \in
\bN_0$.

We consider first the case when $s$ is odd.

Set
$U_{-1}:=(-\infty,0)$, $U_0:=(-1,1)$, $U_1:=(0,\infty)$, that are
invariant sets under $\vp$ and they cover $\bR$. Denote by $T_k$
the restriction of $T=\cf-\lambda I$ to $\cA(U_k)$, $k=-1,0,1$. By
Theorem \ref{spectrum1}, $\ker T_0$ is finite dimensional,
since $\vp$ and $\vp^{[2]}$ have on $U_0$ only one fixed point. On
the other hand, since $\vp$ is a real analytic diffeomorphism both
in $U_{-1}$ and in $U_1$, we can apply Corollary
\ref{cor2spectrum1} to conclude that $\ker T_i$ is finite
dimensional for $i=-1,1$. Now set $U_{-1,0}:=U_{-1} \cap U_0=(-1,0)$ and
$U_{0,1}:=U_0 \cap U_1=(0,1)$ and $T_{i,k}$ for the operator $T$ restricted
to $\mathcal{A}(U_{i,k})$. Each $U_{i,k}$ is invariant under
$\vp$ and $\vp$ is a real analytic diffeomorphism on each $U_{i,k}$
without fixed points. We can apply
Theorem  \ref{eigen}  to obtain that each $\ker T_{i,k}$ is infinite
dimensional.

In the notation of Theorem \ref{theoremBT}, we have that $S: F \rightarrow E$
is not surjective, because the domain is finite
dimensional and the range is infinite dimensional. We can apply Theorem \ref{theoremBT}
to conclude that $\cf-\lambda I: \Ar \to \Ar$ is not surjective (even more is true: $\operatorname{codim}\im \cf-\lambda I=\infty$).

Now suppose that $s$ is even. We define $U_0=\bR\setminus\{0\}$,
$U_1=(-1,1)$. They are invariant with respect to $\vp$. We will
describe $\ker T\vert_{\cA(U_i)}$.

First, if $f\in \ker T\vert_{\cA(U_0)}$, then for $x<0$ we
can calculate $f(x)=(1/\lambda)f(x^s)$ so it is uniquely
determined by $f\vert_{\bR_+}$ and the latter function belongs to
$\ker T\vert_{\cA(\bR_+)}$. Since $\vp$ is a
diffeomorphism and $\vp^{[2]}$ has one fixed point on $\bR_+$ thus
by Corollary \ref{cor2spectrum1}, the kernel is finite
dimensional. Now, $\ker T\vert_{\cA(U_1)}$ is finite
dimensional by Theorem \ref{spectrum1}.

For the space $\cA(U_0\cap U_1)$ the kernel of $T_\lambda$ is
infinite dimensional, since again its elements $f$ are determined
uniquely by $f\vert_{(0,1)}$ and the latter belongs to $\ker
T\vert_{\cA(0,1)}$. By Theorem \ref{eigen}, this
kernel is infinite dimensional.

We apply again Theorem \ref{theoremBT} to conclude that $T=\cf-\lambda I:\Ar\to\Ar$ is not
surjective, since  the map $S:F \rightarrow E$ is not surjective. Indeed,
the  domain is finite
dimensional but the range space is infinite dimensional and so $\operatorname{codim}\im \cf-\lambda I=\infty$.
\qed
\vskip0.5truecm
In order to get more examples, the following observation is useful. The proof is easy.

\begin{lemma}\label{diffchange}
Let $\vp: \bR \to \bR$ be a real analytic self map and let $\delta:\bR \to \bR$ be a
real analytic diffeomorphism. The real analytic self map $\psi:\bR \to \bR, \psi:= \delta^{-1} \circ \vp \circ \delta,$
satisfies $\sigma_p(C_\psi) = \sigma_p(C_\vp)$ and $\sigma(C_\psi) = \sigma(C_\vp)$. Moreover, $u$ is a fixed point
of $\vp$ if and only if $\delta^{-1}(u)$ is a fixed point of $\psi$. If $u$ is a fixed point of $\vp$, then $\psi'(\delta^{-1}(u)) = \vp'(u)$.
Moreover, a subset $A$ of $\bR$ is invariant by $\vp$ if and only if $\delta^{-1}(A)$ is invariant for $\psi$.
\end{lemma}

As an immediate consequence of Lemma \ref{diffchange} and Example \ref{exmonomials}, we get

\begin{example}\rm
Let $\delta: \bR \to \bR$ be a real analytic diffeomorphism and let $s \in \bN$, $s>1$.
The real analytic self map $\vp(x):= \delta^{-1}(\delta(x)^s), x \in \bR,$ satisfies
$\sigma_p(\cf) = \{1 \}$ and $\sigma(\cf)= \bC$.

To mention a concrete explicit example, if $\delta(x)=e^x - e^{-x}$, then $\delta^{-1}(x)= \log(\frac{1}{2}((y^2 + 4)^{1/2} -y))$.
\end{example}

\section{Quadratic function}

We now investigate the spectrum $\sigma(\cf)$ of $\cf: \Ar \to \Ar$, when
$\vp(x)=a x^2 + b x + c, a \neq 0$. We distinguish several cases, depending on the number of fixed
points on $\vp$. Lemma \ref{diffchange} will be used to reduce the problem to the study of certain parameter family. \\

Case 1. The symbol $\vp$ has no fixed points. In this case, since $\vp$ has
only one critical point, it follows from Corollary \ref{fpf} that $\sigma(\cf) = \bC$ and $\sigma_p(\cf)= \bC \setminus \{ 0 \}$. \\

Case 2 and 3. The symbol has fixed points $v$, $u$ (possibly $u=v$) and
$$
\vp(x)=x+a(x-u)(x-v), \qquad a\neq 0.
$$
Choose $u\leq v$ if $a<0$ or $u>v$ if $a>0$. then taking
$$
\delta(x)=-\frac{1}{a}x+u
$$
we will obtain
$$
\psi(x)=\delta^{-1}\circ \vp\circ \delta(x)=\left[ 1+a(u-v)\right] x-x^2.
$$
Clearly, $u$, $v$ were selected to ensure that $\mu:=1+a(u-v)\geq 1$. By Lemma \ref{diffchange}, it is enough to consider
$$
\vp(x)=-x^2+\mu x\qquad \text{with $\mu\geq 1$}.
$$

Case 2. The symbol $\vp$ has only one fixed point $u$. Then the behaviour of $\cf$ is the same as $C_\psi$ for $\psi(x)= -x^2+x$.

\begin{proposition}\label{case2}
In Case 2, $\vp(x)=-x^2+x$, $\sigma_p(\cf)=\{1\}$  but $\sigma(\cf)\supseteq \{0\}\cup[1,+\infty)$.
\end{proposition}

\Proof
Clearly, by Theorem \ref{eigen} and Proposition \ref{spectrum0} we have
$$
\sigma_p(\cf)=\{ 1\}, \qquad \sigma(\cf)\supset \{1,0\}.
$$

Now, we consider the map $T:\Ar\to\Ar$, $T(f)=\cf(f)-\lambda f$, and we will show that  the function $\id$, $\id(x)=x$, does not belong to the image of $T$ for $\lambda\in (1,+\infty)$.

It is a simple calculation that there is the unique formal solution $f(z)=\sum_{n=0}^\infty f_nz^n$ of the equation $T(f)=\id$, with
\begin{equation}\label{fl}
f_0=0, \quad f_1=\frac{1}{1-\lambda}, \quad f_n=\frac{1}{1-\lambda}\sum_{j=1}^{\left[\frac{n}{2}\right]}\left(\begin{matrix} n-j\\ j\end{matrix}\right)(-1)^{j-1}f_{n-j}.
\end{equation}
It is clear that for $1-\lambda<0$, $f_1<0$ and $f_2>0$. Then one shows inductively by (\ref{fl}) that $f_n(-1)^{n-1}>0$. Therefore
$$
\vert f_n\vert>\frac{1}{\vert 1-\lambda\vert}\left(\begin{matrix} n-1\\ 1\end{matrix}\right)\vert f_{n-1}\vert=\frac{n-1}{\vert 1-\lambda\vert}\vert f_{n-1}\vert.
$$
Again inductively
$$
\vert f_n\vert>\frac{(n-1)!}{\vert 1-\lambda\vert^n}
$$
and the series $\sum_{n=0}^\infty f_nz^n$ is nowhere convergent. So there is no real analytic solution $f$ for the equation $T(f)=\id$ for any $\lambda>1$. \qed

\vspace{.2cm}

In Case 2. we  also get  the following partial positive step.

\begin{proposition}\label{formal} Let $\vp(x)=-x^2+x$, $\lambda\neq 0$, and assume that the equation
\begin{equation}\label{eq7}
f(\vp(x))-\lambda f(x)=\gamma(x)
\end{equation} has an analytic solution $f$ on $(-\eps,\eps)$ for some  $\eps>0$ and some fixed $\gamma\in \Ar$.
Then this solution extends to a solution on $\bR$.
\end{proposition}

\begin{remark} Note that (\ref{eq7}) has always the unique formal solution --- we do not know for precisely which $\lambda$ and $\gamma$ it converges around zero.
\end{remark}

\Proof We observe that for $x\in (0,1)$ $\vp^{[n]}(x)\to 0$ as $n\to \infty$ so, by Lemma \ref{ext}, we extend $f$ to $(-\eps, 1)$.

Since $\vp:(-\infty, 1/2)\to (-\infty, 1/4)$ is a diffeomorphism
onto, its inverse $\psi$ is a diffeomorphism $(-\infty,
1/4)\mapsto (-\infty, 1/2)$. Moreover, for $x\in (-\infty, 1/2)$
the equation
$$
f(\vp(x))-\lambda f(x)=\gamma(x)
$$
is equivalent to
$$
f(\psi(x))-(1/\lambda)f(x)=-(1/\lambda)\gamma(\psi(x)).
$$
Since for $x\in (-\infty, 0)$ $\psi^{[n]}(x)\to 0$ as $n\to
\infty$ we can extend $f$ onto $(-\infty, 0)$ by Lemma \ref{ext}.

Summarizing, $f$ is a solution of (\ref{eq7}) on $(-\infty, 1)$.

Since for $x\in (0, 1/2)$ and $y=1-x\in (1/2, 1)$ we have
$$
\lambda f(x)+\gamma(x)=f(\vp(x))=f(\vp(y))=\lambda f(y)+\gamma(y)
$$
we have
$$
f(y)=f(1-y)+\frac{\gamma(1-y)-\gamma(y)}{\lambda}.
$$
This formula extends $f$ analytically onto $(1/2, +\infty)$ and so
on the whole real line.  Clearly, this extension satisfies
(\ref{eq7}) everywhere. \qed

\vspace{0.2cm}

Case 3. The symbol $\vp$ has exactly two different fixed points $u$, $v$. Then the behaviour of $\cf$ is the same as $C_\psi$ for
$$
\psi(x)=-x^2+\mu x,\qquad \mu>1.
$$
The map $\psi$ has one critical point $\mu/2$ and fixed points $0$ and $\mu-1$.

Case 3.1. Critical point is outside the interval joining fixed points, i.e., $\mu/2>\mu-1$ so $1<\mu<2$.

Case 3.2. Critical point is equal to one fixed point, i.e., $\mu/2=\mu-1$ so $\mu=2$.

Case 3.3. Critical point is between fixed points, i.e., $\mu/2<\mu-1$ so $\mu>2$.

\begin{proposition} In Case 3.1 and 3.2, $\vp(x)=-x^2+\mu x$, $1<\mu\leq 2$, $\sigma(\cf)=\bC$ but
$\sigma_p(\cf)=\{ 1\}$ and for every $\lambda \in \bC$ kernel of $\cf-\lambda I$ is one-dimensional.
\end{proposition}

\Proof We define
$$
U_1:=(-\infty, \mu-1)\cup (1,+\infty), \qquad U_0:= (0,\mu).
$$
Observe that for every $x\in (-\infty, \mu-1)$ there is exactly
one point $y\in (1, \infty)$ such that $\vp(x)=\vp(y)$
($y=\mu-x$). Thus $f\in \cA(U_1)$ belongs to the kernel of
$\cf-\lambda I$, $\lambda\neq 0$ if and only if $f\vert_{(-\infty,
\mu-1)}\in \ker (\cf-\lambda I)$ on $\cA(-\infty, \mu-1)$ and
$$
f(y)=f(\mu-y),\qquad y\in (1,+\infty).
$$
Therefore the kernels of $\cf -\lambda I$ on $\cA(U_1)$ and on
$\cA(-\infty, \mu-1)$ have the same dimension. In the latter case
$\vp$ and $\vp^{[2]}$ have only one fixed point $0$ on $(-\infty,
\mu-1)$. Hence, by Theorem \ref{eigen}, $\ker \cf-\lambda I$ is
one dimensional in $\cA(U_1)$.

Since $\vp:(0,\mu)\to (0,\mu)$ and its square have exactly one
fixed point $\mu-1$,   $\ker \cf-\lambda I$ is one dimensional
in $\cA(U_0)$ by Theorem \ref{eigen}.

Finally, $U_1\cap U_0=(0,\mu-1)\cup (1,\mu)$.  Again as in case of
$U_1$ one can prove that the kernels of $\cf-\lambda I$ in
$\cA(U_1\cap U_0)$ and in $\cA(0,\mu-1)$ are isomorphic. Since
$\vp$ has no fixed point in $(0, \mu-1)$ and only one critical
point, the required kernel is
infinite dimensional by Theorem \ref{eigen}.

We apply Theorem \ref{theoremBT} to conclude that,  for every $\lambda\neq 0$, the map
$\cf-\lambda I:\Ar\to \Ar$ is not surjective. By Proposition
\ref{spectrum0}, $\sigma(\cf)=\bC$.

The remaining parts of the Proposition follows from Theorem \ref{eigen}. \qed

\vspace{.2cm}

In Case 3.3 we have only the following partial result:

\begin{proposition} If $\vp(x)=-x^2+\mu x$, $\mu> 2$, then $\sigma(\cf)\supset \{\lambda\in \bC:\vert \lambda \vert \leq1\} \cup \{ \mu^n \}_n \cup \{ (2-\mu)^n \}_n$.
\end{proposition}

\Proof Clearly,  $\{ 0, 1 \} \cup \{ \mu^n \}_n \cup \{ (2-\mu)^n \}_n \subset \sigma(\cf)$ by Propositions \ref{spectrum0} and \ref{eigen}.

Fix $\lambda \neq 1$, $0<\vert \lambda\vert \leq 1$. Since $\vp(1)=\mu - 1$, we can select a sequence  $x_n\in (0, 1) \subset (0,\mu)$, $x_1=1$, such that
$$
\vp^{[n]}(x_n)=\mu-1, \ \ n \in \bN
$$
Such $x_n$ exists because $\vp(\mu/2)=\mu^2/4\geq \mu$ and $1 < \mu/2$. In fact,
$x_n=\frac{\mu-\sqrt{\mu^2-4x_{n-1}}}{2} = \frac{2 x_{n-1}}{\mu+\sqrt{\mu^2-4x_{n-1}}}, n \in \bN$.
Hence $x_n/x_{n-1} < 2/\mu < 1, n \in \bN,$ and  $(x_n)_\nN$ is decreasing
to zero as $n\to \infty$. Moreover, $x_n/x_{n-1}$
tends to $1/\mu$ as $n\to \infty$.

If $\cf f- \lambda f= \gamma$ for some $\gamma\in \Ar$, then
$$
\lambda f(\mu-1)+\gamma(\mu-1)=f(\mu-1), \quad \text{so}\quad
f(\mu-1)=\gamma(\mu-1)/(1-\lambda),
$$
and analogously
$$
f(0)=\frac{\gamma(0)}{1-\lambda},
$$
because both $0$ and $\mu-1$ are fixed points for $\vp$.  Since
$$
f\circ \vp^{[n]}(x)=\lambda^nf(x)+\lambda^{n-1}\gamma(x)+\lambda^{n-2}\gamma(\vp(x))+\dots+\gamma(\vp^{[n-1]}(x)),
$$
we have
\begin{equation}\label{eq6}
\frac{\gamma(\mu-1)}{1-\lambda}=\lambda^{n}f(x_n)+\lambda^{n-1}\gamma(x_n)+\dots+\lambda\gamma(\vp^{[n-2]}(x_n))+\gamma(1).
\end{equation}
Fix $c<1/6$. If $\gamma(0)=0$, we can choose $n_1 > 2$ so that, for $n \geq n_1$,  $\vert
x_n\vert \leq c$ and $\vert
f(x_n)\vert \leq c$, since $f(x_n)$ tends to
$f(0)=\gamma(0)/(1-\lambda)=0$.

Take $\gamma(x)=x^k\gamma_0(x)$, with $\gamma_0 \in \Ar$ satisfying $\gamma_0(\mu-1)=0$
and $1/2<\vert \gamma_0(x)\vert <1$ for $x\in [0,1)$,
$\gamma_0(1)=1$. Since $x_n/x_{n-1}$
tends to $1/\mu$ as $n\to \infty$, there is $n_0>n_1$ such that
$x_n<(1/2)x_{n-1}$ for $n>n_0$. We can take $k$ big enough to ensure that $\vert
\gamma(x_j)\vert<c/n_0$ for $2 \leq j \leq n_0$. Then, for $n>n_0$, we have
$\vert\gamma(x_n)\vert\leq 2c(1/2)^{n-n_0}$. Summarizing:
$$
\vert\lambda^{n}f(x_n)+\lambda^{n-1}\gamma(x_n)+\dots+\lambda\gamma(\vp^{[n-2]}(x_n))\vert\leq
2c+2c\frac{1}{1-(1/2)^k}\leq 6c.
$$
As $c<1/6$, we reach a contradiction with (\ref{eq6}). Therefore
$\cf-\lambda I$ is not surjective. \qed

\vspace{.5cm}

{\bf Acknowledgement:}  The research of the authors was partially
supported by MEC and FEDER Project MTM2013-43540-P and the work of
of Bonet by the grant GV Project Prometeo II/2013/013. The research of Doma{\'n}ski
was supported   by National
Center of Science, Poland, grant no. DEC-2013/10/A/ST1/00091.

%}

%{\bf Authors' Addresses:}
%
%\vspace{2mm}
%
%\begin{minipage}{7.5cm}
%J.\ Bonet
%
%Instituto Universitario de Matem\'{a}tica Pura y Aplicada IUMPA
%
%Universitat Polit\`{e}cnica de Val\`{e}ncia
%
%E-46071 Valencia, SPAIN
%
%e-mail: jbonet@mat.upv.es
%\end{minipage}
%\hspace*{\fill}
%\parbox{7.5cm}{P. Doma{\'n}ski
%
%Faculty of Mathematics and Comp. Sci.
%
%A. Mickiewicz University Pozna{\'n}
%
%Umultowska 87
%
%61-614 Pozna{\'n}, POLAND
%
%e-mail: domanski@amu.edu.pl
%
%}

\end{document}